\documentclass[11pt]{article}
\usepackage[a4paper]{anysize}\marginsize{3.5cm}{3.5cm}{1.3cm}{2cm}
\pdfpagewidth=\paperwidth \pdfpageheight=\paperheight
\usepackage{amsfonts,amssymb,amsthm,amsmath,eucal}
\usepackage{pgf}
\usepackage{bbm,array}
\usepackage{tikz} 
\usepackage{float}
\restylefloat{table}
\usepackage{subfigure}
\usepackage{caption}
\usetikzlibrary{arrows}
\usepackage{multicol}
\usepackage{multirow}

\theoremstyle{plain}
\newtheorem{thm}{Theorem}[section]

\theoremstyle{definition}

\newtheorem{thevarthm}[thm]{\varthmname}

\newenvironment{varthm*}[1]{\trivlist\item[]{\bf #1.}\it}{\endtrivlist}

\newcommand\be{\begin{eqnarray*}}
\newcommand\ee{\end{eqnarray*}}

\newcommand\C{\mathbb C}

\renewcommand\P{\mathbb P}

\def\field{\C}
\newcommand\newop[2]{\def#1{\mathop{\rm #2}\nolimits}}
\newop\mod{mod}
\newop\log{log}
\newop\ord{ord}
\newop\Gal{Gal}
\newop\SL{SL}
\newop\GL{GL}
\newop\Bl{Bl}
\newop\mult{mult}
\newop\mass{mass}
\newop\div{div}
\newop\codim{codim}
\newop\sing{sing}
\newop\vdim{vdim}
\newop\edim{edim}
\newop\Ass{Ass}
\newop\size{size}
\newop\reg{reg}
\newop\areg{areg}
\newop\asreg{asreg}
\newop\satdeg{satdeg}
\newop\supp{supp}
\newop\gin{gin}
\newop\ini{in}
\newop\vol{vol}
\newop\sat{sat}
\newop\length{length}
\newop\depth{depth}
\newop\characteristic{char}

\def\keywordname{{\bfseries Keywords}}%
\def\keywords#1{\par\addvspace\medskipamount{\rightskip=0pt plus1cm
\def\and{\ifhmode\unskip\nobreak\fi\ $\cdot$
}\noindent\keywordname\enspace\ignorespaces#1\par}}
\def\subclassname{{\bfseries Mathematics Subject Classification
(2000)}\enspace}
\def\subclass#1{\par\addvspace\medskipamount{\rightskip=0pt plus1cm
\def\and{\ifhmode\unskip\nobreak\fi\ $\cdot$
}\noindent\subclassname\ignorespaces#1\par}}

\definecolor{ttqqqq}{rgb}{0.,0.,0.}
\definecolor{zzttqq}{gray}{0.4}

\begin{document}

\author{\L .~Farnik,  J.~Kabat, M.~Lampa-Baczy\'nska, H.~Tutaj-Gasi\'nska}
\title{Containment problem and combinatorics}
\date{\today}
\maketitle
\thispagestyle{empty}

\begin{abstract}
In this note we consider two configurations of twelve lines with nineteen
triple points (i.e., points where three lines meet). Both of them have the same combinatorial features. In both
configurations nine of twelve lines have five triple points and one double point,
and the remaining three lines have four triple points and three double points.
Taking the ideal of the triple points of these
configurations we discover that, quite surprisingly, for one of the
configurations the containment $I^{(3)} \subset I^2$ holds, while for the
other it does not. Hence for ideals of points defined by
configurations of lines the (non)containment of a
symbolic power in an ordinary power is not determined alone by
combinatorial features of the arrangement.

Moreover, for the configuration with the non-containment $I^{(3)} \nsubseteq I^2$
we examine its parameter space, which turns out to be a rational curve,
and thus establish the existence of a rational non-containment configuration
of points. Such rational examples are very rare.

 \subclass{52C30,  14N20, 05B30}
\end{abstract}


\section{Introduction}\label{intro}

The notion of the symbolic power of an ideal appears recently in many problems. Let  $I\subset \field[\mathbb{P}^N]=\field[x_0,\dots,x_N]$  be a homogeneous
ideal. By $m$-th symbolic power of $I$ we mean $ I^{(m)} = \field[\mathbb{P}^N] \cap \left( \bigcap_{\mathfrak{p} \in \Ass(I)} (I^{m})_{\mathfrak{p}} \right). $
For
a radical ideal  $I$, the Nagata-Zariski theorem says that $I^{(m)}$ is the ideal of all $f\in I$ which vanish to order at least
$m$ along the zero-set of $I$. The main question concerning the symbolic powers may be stated as follows.
 For which $r$ and $m$ does the containment
$$
I^{(m)}\subset I^r
$$
hold? Or, more generally, when $M=(x_0,\ldots,x_N)$, for which $r$, $m$ and $j$ do we have
$$
I^{(m)} \subset M^jI^r?
$$
Ein, Lazarsfeld and Smith \cite{ELS},  and Hochster with Huneke \cite{HoHu} showed that, for any ideal $I \subset \field[\mathbb{P}^N]$, the containment $I^{(rN)} \subset
I^r$ holds.

 The following question was asked by Harbourne and Huneke a few years ago. Let
$I$ be an ideal of points in $\mathbb{P}^N$. Does then the containment
$$
I^{(rN-(N-1))} \subset I^r
$$
hold for all $r$? Lots of examples suggested that the answer is positive.
 For an ideal of points in $\mathbb{P}^2$ in particular, the question was if $I^{(2r-1)} \subset I^r$ holds.  In the paper \cite{DST13}, the first
counterexample for the case $r=2$ was presented. Since then quite a few counterexamples appeared, see, e.g., \cite{CGMLLPS2015}, \cite{HarSec13}, \cite{SS},  \cite{MJ2015},
 \cite{malara-szpond}, \cite{MS} or are announced \cite{BFFHTG}.
The case $r>2$ is still open.

The first real -- and rational -- counterexamples (i.e., counterexamples where the coordinates of all points are real numbers)
 come from \cite{CGMLLPS2015}, \cite{res} and  \cite{MJ2015}. They are modifications of  B\"or\"oczky configuration of $12$ lines. The non-existence of a rational counterexample among  B\"or\"oczky configurations of $13$,  $14$, $16$, $18$ and $24$ lines is studied in \cite{FKL-BT-G}.
Recently a new rational counterexample appeared, see \cite{malara-szpond}.

In the paper of Bokowski and Pokora, \cite{B-P}, two non-isomorphic (and non-isomor\-phic to B\"or\"oczky)   examples of real configurations are considered. They are named there $C_2$ and $C_7$.

In this paper we consider the two configurations,  $C_2$ and $C_7$. These configurations  are realizable over the reals, and, what is interesting, they have the
same combinatorial features as  B\"or\"oczky configuration of $12$ lines, i.e., the lines intersect in $19$ triple points, $9$ lines have $5$ triple points on
them and $3$ lines have $4$ triple points. We describe the parameter spaces of these configurations. It turns out, that one of them, $C_2$,  is ``rigid", this
means that fixing some four out of $19$ triple points (by a projective automorphism) to be $(1:0:0), (0:1:0), (0:0:1), (1:1:1)$, the coordinates of other
points can be computed, and these coordinates are non-rational. Moreover, for this configuration the containment $I_2^{(3)}\subset I_2^2$ holds, where $I_2$ is
the radical ideal of the triple points of the configuration. The second configuration, namely $C_7$, turns out to have a one dimensional projective space as a
parameter space. Thus, we can take all the triple points of the configuration with rational coefficients. The radical ideal of these points, $I_7$ gives a new rational example of the non-containment $I_7^{(3)}\nsubseteq I_7^2$.

\section{Configuration $C_2$}\label{C2}

The real realization of the  configuration $C_2$ is pictured below. Points $P_1$ and $P_2$ are ``at infinity".

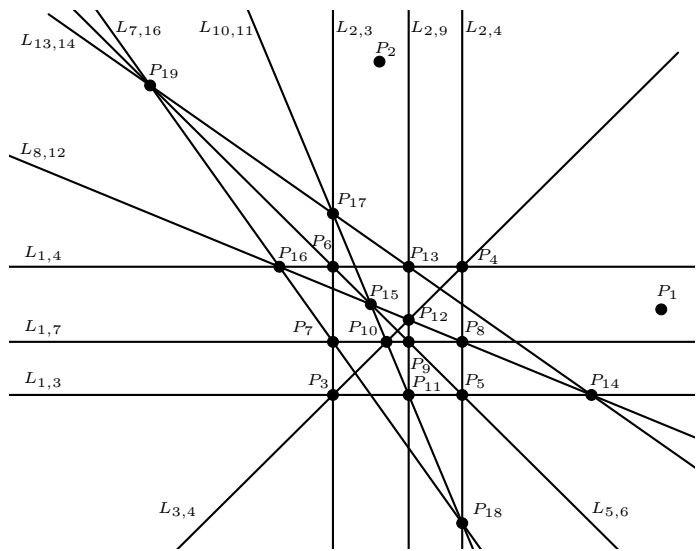
\begin{figure}[h]
\centering
\begin{tikzpicture}[line cap=round,line join=round,x=1.0cm,y=1.0cm,scale=1.7]
\clip(-2.5,-1.2) rectangle (2.8,3.);
\draw [line width=0.8pt,domain=-2.6142181668765607:7.7791046265004] plot(\x,{(-1.-0.*\x)/-1.});
\draw [line width=0.8pt,domain=-2.6142181668765607:7.7791046265004] plot(\x,{(--0.17157287525381-0.*\x)/0.41421356237309515});
\draw [line width=0.8pt,domain=-2.6142181668765607:7.7791046265004] plot(\x,{(-0.-0.*\x)/1.});
\draw [line width=0.8pt] (0.,-2.023398296524389) -- (0.,3.17784973864309);
\draw [line width=0.8pt] (0.5857864376269049,-2.023398296524389) -- (0.5857864376269049,3.17784973864309);
\draw [line width=0.8pt] (1.,-2.023398296524389) -- (1.,3.17784973864309);
\draw [line width=0.8pt,domain=-2.6142181668765607:7.7791046265004] plot(\x,{(--1.-1.*\x)/1.});
\draw [line width=0.8pt,domain=-2.:2.67] plot(\x,{(-0.--1.*\x)/1.});
\draw [line width=0.8pt,domain=-2.6142181668765607:7.7791046265004] plot(\x,{(-0.3431457505076197--0.1715728752538097*\x)/-0.41421356237309515});
\draw [line width=0.8pt,domain=-2.6142181668765607:7.7791046265004] plot(\x,{(--0.24264068711928516-0.41421356237309515*\x)/0.1715728752538097});
\draw [line width=0.8pt,domain=-2.2:7.7791046265004] plot(\x,{(--2.-1.*\x)/1.414213562373097});
\draw [line width=0.8pt,domain=-2.6142181668765607:7.7791046265004] plot(\x,{(--0.17157287525381062-0.5857864376269049*\x)/0.4142135623730967});
\begin{tiny}
\draw [fill=black] (1.,1.) circle (1.2pt);
\draw[color=black] (1.2,1.08) node {$P_{4}$};
\draw [fill=black] (1.,0.) circle (1.2pt);
\draw[color=black] (1.1,0.11) node {$P_{5}$};
\draw [fill=black] (0.,0.) circle (1.2pt);
\draw[color=black] (-0.11,0.11) node {$P_{3}$};
\draw [fill=black] (0.,1.) circle (1.2pt);
\draw[color=black] (-0.085,1.19) node {$P_{6}$};
\draw [fill=black] (1.,0.41421356237309515) circle (1.2pt);
\draw[color=black] (1.1,0.52) node {$P_{8}$};
\draw [fill=black] (0.5857864376269049,0.41421356237309515) circle (1.2pt);
\draw[color=black] (0.68,0.23) node {$P_{9}$};
\draw [fill=black] (0.5857864376269049,0.) circle (1.2pt);
\draw[color=black] (0.73,0.09) node {$P_{11}$};
\draw[color=black] (-2.24,1.1) node {$L_{1,4}$};
\draw[color=black] (-2.24,0.52) node {$L_{1,7}$};
\draw[color=black] (-2.24,0.11) node {$L_{1,3}$};
\draw[color=black] (0.17,2.85) node {$L_{2,3}$};
\draw[color=black] (0.75,2.85) node {$L_{2,9}$};
\draw[color=black] (1.17,2.85) node {$L_{2,4}$};
\draw [fill=black] (0.5857864376269049,1.) circle (1.2pt);
\draw[color=black] (0.71,1.1) node {$P_{13}$};
\draw [fill=black] (0.,0.41421356237309515) circle (1.2pt);
\draw[color=black] (-0.23,0.52) node {$P_{7}$};
\draw[color=black] (2.15,-0.9) node {$L_{5,6}$};
\draw[color=black] (-1.2,-0.9) node {$L_{3,4}$};
\draw [fill=black] (0.5857864376269049,0.5857864376269049) circle (1.2pt);
\draw[color=black] (0.78,0.63) node {$P_{12}$};
\draw [fill=black] (0.41421356237309515,0.41421356237309515) circle (1.2pt);
\draw[color=black] (0.235,0.52) node {$P_{10}$};
\draw[color=black] (-2.24,1.9) node {$L_{8,12}$};
\draw[color=black] (-0.82,2.85) node {$L_{10,11}$};
\draw [fill=black] (0.292893218813453,0.707106781186547) circle (1.2pt);
\draw[color=black] (0.3992033138157002,0.8) node {$P_{15}$};
\draw [fill=black] (2.,0.) circle (1.2pt);
\draw[color=black] (2.1054328279976198,0.11) node {$P_{14}$};
\draw [fill=black] (-0.4142135623730967,1.) circle (1.2pt);
\draw[color=black] (-0.3071390119477826,1.1) node {$P_{16}$};
\draw [fill=black] (0.,1.4142135623730967) circle (1.2pt);
\draw[color=black] (0.15,1.52) node {$P_{17}$};
\draw [fill=black] (1.,-1.) circle (1.2pt);
\draw[color=black] (1.2,-0.9) node {$P_{18}$};
\draw[color=black] (-2.2,2.75) node {$L_{13,14}$};
\draw[color=black] (-1.5,2.85) node {$L_{7,16}$};
\draw [fill=black] (-1.4142135623730885,2.4142135623730883) circle (1.2pt);
\draw[color=black] (-1.3070262003662192,2.53) node {$P_{19}$};
\draw [fill=black] (2.54,0.6674296111460296) circle (1.2pt);
\draw[color=black] (2.59,0.8) node {$P_1$};
\draw [fill=black] (0.36,2.6) circle (1.2pt);
\draw[color=black] (0.41,2.7) node {$P_2$};
\end{tiny}
\end{tikzpicture}
\end{figure}

By a projective automorphism we may move any four general points of $\mathbb{P}^2$ into other four general points. Thus, we may assume (with the
notation as in the picture) that $P_{1}=(1:0:0), P_{2}=(0:1:0), P_3=(0:0:1)$ and $P_4=(1:1:1)$. We take the following lines:
$$
\begin{array}{ll}
L_{1,3}: \ y=0,\\
L_{2,4}: \ x-z=0,\\
L_{1,4}: \ y-z=0,\\
L_{3,4}: \ x-y=0,\\
L_{2,3}: \ x=0,
\end{array}
$$
where $L_{i,j}$, is the line passing through the points $P_{i}$ and $P_{j}$. Then we obtain the points
$$
\begin{array}{ll}
P_5= L_{1,3} \cap L_{2,4}=(1,0,1),\\
P_6= L_{1,4} \cap L_{2,3}=(0,1,1)
\end{array}
$$
and the line
$$L_{5,6}: \ x+y-z=0.$$
We need now to introduce a parameter to proceed with the construction. Thus we take the point $P_7 = (0,1,a) \in L_{2,3}$. Since all points and lines in the configuration should be distinct, we assume that $a \neq 1$ and $a \neq 0$. We obtain the remaining lines and points in the following order:
$$L_{1,7}:\ z-ay=0,$$
$$
\begin{array}{ll}
P_8= L_{1,7} \cap L_{2,4}=(a,1,a),\\
P_9= L_{1,7} \cap L_{5,6}=(a-1,1,a),\\
P_{10}= L_{1,7} \cap L_{3,4}=(1,1,a),\\
\end{array}
$$
$$L_{2,9}:\ ax-(a-1)z=0, $$
$$
\begin{array}{ll}
P_{11}= L_{2,9} \cap L_{1,3}=(a-1,0,a),\\
P_{12}= L_{2,9} \cap L_{3,4}=(a-1,a-1,a),\\
P_{13}= L_{2,9} \cap L_{1,4}=(a-1,a,a),
\end{array}
$$
$$L_{8,12}:\ a(2-a)x-ay+(a-1)^2z=0, $$
$$
\begin{array}{ll}
P_{14}= L_{8,12} \cap L_{1,3}=((a-1)^2,0,a(a-2)),\\
P_{15}= L_{8,12} \cap L_{5,6}=(a^2 -3a+1,-1,a(a-3)),\\
P_{16}= L_{8,12} \cap L_{1,4}=(a^2-3a+1,a(a-2),a(a-2)),
\end{array}
$$
$$L_{10,11}:\ ax-a(a-2)y-(a-1)z=0, $$
$$
\begin{array}{ll}
P_{17}= L_{10,11} \cap L_{2,3}=(0,a-1,a(a-2)),\\
P_{18}= L_{10,11} \cap L_{2,4}=(a(2-a),1,a(2-a)),
\end{array}
$$
$$L_{13,14}:\ a(a-2)x+(a-1)y-(a-1)^2 z=0, $$
$$L_{7,16}:\ a(a-1)(a-2)x-a(a^2-3a+1)y+ (a^2-3a +1) z=0, $$
$$P_{19}= L_{13,14} \cap L_{7,16}=(a^5 -5a^4 +7 a^3 -a^2 -3a +1, a^3 (a-2)^2, a^5 -4a^4 +3a^3 +3a^2 -2a).$$
Almost all points in the configuration are triple directly from the construction.  Only in four of them, i.e., $P_{15}$, $P_{17}$, $P_{18}$ and $P_{19}$, we must verify this fact. We need to check the following incidences:
$$
\begin{array}{ll}
P_{15}= L_{8,12} \cap L_{5,6} \cap L_{10,11},\\
P_{17}= L_{10,11} \cap L_{2,3}  \cap L_{13,14},\\
P_{18}= L_{10,11} \cap L_{2,4} \cap L_{7,16},\\
P_{19}= L_{13,14} \cap L_{7,16} \cap L_{5,6}.
\end{array}
$$
By the determinant condition we conclude that the lines  $L_{8,12}$, $L_{5,6}$ and $L_{10,11}$ always meet at a point, but the remaining incidences occur under
the algebraic condition
$$a^2 -2a -1=0.$$
Thus,  the configuration has no rational realization.

Then, implementing, e.g. in Singular (\cite{DGPS}), the ideal $I_2$ of all the triple points we check that $I_2^{(3)}\subset I_2^2$. This inclusion may be explained also more
theoretically. From  \cite{BocHar10a}, we have that if $\alpha(I^{(m)})\geq r\cdot  \textrm{reg}I$ (where $\alpha(J)$ denotes the least degree of a nonzero form in
a homogeneous ideal $J$), then the containment $I^{(m)}\subset I^r$ holds. It may be computed (e.g., with Singular) that reg $I_2=6$ and $\alpha(I^{(3)})=12$.
Thus,  $I_2^{(3)}\subset I_2^2$.

There is an interesting phenomenon  that for ideal $I_{2}$ the inclusion $I_2^{(3)}\subset I_2^2$ is true, while for other configurations of $12$  lines,
B\"{o}r\"{o}czky and $C_7$, with  the same combinatorial features, the inclusion does not occur, see  the next section for $C_7$ and  \cite{MJ2015} for
B\"{o}r\"{o}czky. Thus, the combinatoric features of the arrangement do not determine the containment.

\section{Configuration $C_7$}\label{C7}

The real realization of the  configuration $C_7$ is in the picture (the points $P_1, P_2, P_3$ are ``at infinity"):

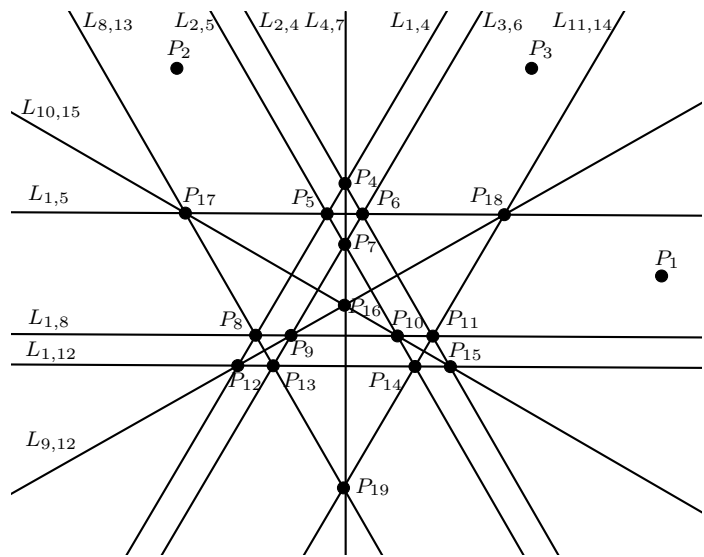
\begin{figure}[h]
\centering
\begin{tikzpicture}[line cap=round,line join=round,x=1.0cm,y=1.0cm,scale=0.15]
\clip(-31.5,-21) rectangle (29,27);
\draw [line width=0.8pt,domain=-45.90951793747096:30.415852243120977] plot(\x,{(-48.66002470361438-8.051009098402481*\x)/-4.701789313467024});
\draw [line width=0.8pt,domain=-45.90951793747096:30.415852243120977] plot(\x,{(-15.359111643461734-0.04636443950339797*\x)/9.32327306204971});
\draw [line width=0.8pt,domain=-45.90951793747096:30.415852243120977] plot(\x,{(-11.260615728558832--8.097373537905879*\x)/-4.621483748582687});
\draw [line width=0.8pt,domain=-45.90951793747096:30.415852243120977] plot(\x,{(-19.226291119539457-4.621483748582687*\x)/-8.097373537905879});
\draw [line width=0.8pt,domain=-45.90951793747096:30.415852243120977] plot(\x,{(-2.3662677197417743-4.701789313467024*\x)/8.051009098402481});
\draw [line width=0.8pt] (-2.2,-23) -- (-2.2,27);
\draw [line width=0.8pt,domain=-45.90951793747096:30.415852243120977] plot(\x,{(-3045.243807478152-3.4903035109415743*\x)/701.853684644548});
\draw [line width=0.8pt,domain=-45.90951793747096:30.415852243120977] plot(\x,{(--28.337963708017178-0.015454813167798065*\x)/3.1077576873499035});
\draw [line width=0.8pt,domain=-45.90951793747096:30.415852243120977] plot(\x,{(--17.237659354632964-2.683669699467494*\x)/-1.5672631044890086});
\draw [line width=0.8pt,domain=-45.90951793747096:30.415852243120977] plot(\x,{(-24.5-2.683669699467494*\x)/-1.5672631044890086});
\draw [line width=0.8pt,domain=-45.90951793747096:30.415852243120977] plot(\x,{(-29.7041290129848-2.6991245126352927*\x)/1.5404945828608945});
\draw [line width=0.8pt,domain=-45.90951793747096:30.415852243120977] plot(\x,{(12-2.6991245126352927*\x)/1.5404945828608945});
\begin{scriptsize}
\draw [fill=black] (-2.2839670155268523,6.4383543078910535) circle (15pt);
\draw[color=black] (-0.5,6.5) node {$P_7$};
\draw [fill=black] (-6.985756328993876,-1.6126547905114283) circle (15pt);
\draw[color=black] (-6.,-2.8) node {$P_9$};
\draw [fill=black] (2.337516733055834,-1.6590192300148263) circle (15pt);
\draw[color=black] (3.24,-0.2) node {$P_{10}$};
\draw[color=black] (11.6,26) node {$L_{3,6}$};
\draw[color=black] (-28.2,-0.2) node {$L_{1,8}$};
\draw[color=black] (-15.4,26) node {$L_{2,5}$};
\draw[color=black] (-28,-11) node {$L_{9,12}$};
\draw[color=black] (-28,18.3) node {$L_{10,15}$};
\draw[color=black] (-4,26) node {$L_{4,7}$};
\draw[color=black] (-8,26) node {$L_{2,4}$};
\draw[color=black] (-28,-3.1) node {$L_{1,12}$};
\draw[color=black] (3.5,26) node {$L_{1,4}$};
\draw [fill=black] (-11.660777120832789,-4.280869676811123) circle (15pt);
\draw[color=black] (-11.,-5.7) node {$P_{12}$};
\draw [fill=black] (-8.553019433482882,-4.296324489978922) circle (15pt);
\draw[color=black] (-6.2,-5.7) node {$P_{13}$};
\draw [fill=black] (-10.093514016343777,-1.597199977343629) circle (15pt);
\draw[color=black] (-12.2,-0.3) node {$P_8$};
\draw [fill=black] (3.8780113159167304,-4.358143742650119) circle (15pt);
\draw[color=black] (1.3,-5.65) node {$P_{14}$};
\draw [fill=black] (6.985769003266634,-4.373598555817918) circle (15pt);
\draw[color=black] (8.3,-3) node {$P_{15}$};
\draw [fill=black] (5.445274420405739,-1.6744740431826257) circle (15pt);
\draw[color=black] (8.1,-0.2) node {$P_{11}$};
\draw [fill=black] (-0.7167039110378449,9.122024007358547) circle (15pt);
\draw[color=black] (1.6,10.4) node {$P_6$};
\draw [fill=black] (-3.8244615983877486,9.137478820526345) circle (15pt);
\draw[color=black] (-5.9,10.5) node {$P_5$};
\draw [fill=black] (-2.2571984938987395,11.821148519993839) circle (15pt);
\draw[color=black] (-0.4,12.3) node {$P_4$};
\draw[color=black] (-28.2,10.7) node {$L_{1,5}$};
\draw[color=black] (18.6,26) node {$L_{11,14}$};
\draw[color=black] (-23,26) node {$L_{8,13}$};
\draw [fill=black] (-2.391041102039304,-15.092822540520094) circle (15pt);
\draw[color=black] (0.2,-15) node {$P_{19}$};
\draw [fill=black] (-16.255492347787353,9.199298073197538) circle (15pt);
\draw[color=black] (-15,10.7) node {$P_{17}$};
\draw [fill=black] (11.714326838361776,9.060204754687355) circle (15pt);
\draw[color=black] (10.2,10.7) node {$P_{18}$};
\draw [fill=black] (-2.3107355371549656,1.055560095788266) circle (15pt);
\draw[color=black] (-0.5,0.9) node {$P_{16}$};
\draw [fill=black] (25.5,3.641076073060212) circle (15pt);
\draw[color=black] (26,5.1) node {$P_1$};
\draw [fill=black] (-17.,22) circle (15pt);
\draw[color=black] (-16.8,23.703744806244465) node {$P_2$};
\draw [fill=black] (14.104030273097319,22) circle (15pt);
\draw[color=black] (14.932705720772327,23.703744806244465) node {$P_3$};
\end{scriptsize}
\end{tikzpicture}
\end{figure}

Here, using  a projective automorphism, we may assume (with the notation as in the picture) that $P_1=(1,0,0), P_2=(-1,1,0), P_3=(1,1,0)$
 and $P_4=(0,0,1)$. Then we have lines:
$$
\begin{array}{ll}
L_{2,4}:\ x+y=0,\\
L_{1,4}:\ x-y=0.
\end{array}
$$
We need now to introduce the parameter to proceed with the construction, so take a point on the line $L_{1,4}$:
$$P_5=(a,a,1),$$
where $a\neq 0$. We get the lines
$$
\begin{array}{ll}
L_{1,5}:\ y-az=0,\\
L_{2,5}:\ x+y-2az=0
\end{array}
$$
and the point
$$P_6=L_{2,4}\cap L_{1,5}=(-a,a,1)$$
and then the line
 $$L_{3,6}:\ x-y+2az=0.$$
 To continue we need to choose next point. We take a point on the line $L_{2,5}$.
 $$P_7=(b,-b+2a,1).$$
 We get the line
 $$L_{4,7}:\ 2ax-bx-by=0.$$
 The condition for the lines $L_{4,7}, L_{2,5}, L_{3,6}$ to meet at $P_7$ is
 $$ba=0.$$
 As $a\neq 0$, we have to take $b=0$. Thus, from now on:
 $$P_7=(0,2a,1)$$
and
 $$L_{4,7}:\ 2ax=0.$$
 Again, we need a new parameter. Take a point on the line $L_{1,4}$
 $$P_8=(c,c,1),$$
 where $a\neq c, c\neq 0$.
 Then
 $$L_{1,8}:\ y-cz=0,$$
$$
\begin{array}{l}
P_9=L_{1,8}\cap L_{3,6}=(-2a+c,c,1),\\
P_{10}=L_{1,8}\cap L_{2,5}=(2a-c,c,1),\\
P_{11}=L_{1,8}\cap L_{2,4}=(-c,c,1).
\end{array}
$$

 Now, choose the last parameter by taking a point, again on the line $L_{1,4}$:
 $$P_{12}=(d,d,1),$$
with $d$ different from $0, a$ and $ c$. Then
 $$L_{1,12}:\ y-dz=0,$$
$$
\begin{array}{l}
$$P_{13}=L_{1,12}\cap L_{3,6}=(-2a+d,d,1),\\
$$P_{14}=L_{1,12}\cap L_{2,5}=(2a-d,d,1),\\
$$P_{15}=L_{1,12}\cap L_{2,4}=(-d,d,1),\\
$$L_{10,15}: (c-d)x+(c-d-2a)y+2adz=0,\\
$$P_{17}=L_{10,15}\cap L_{1,5}=(2a^2-ac-ad,ac-ad,c-d),\\
$$L_{9,12}: (c-d)x+(2a-c+d)y-2adz=0,\\
$$P_{18}=L_{9,12}\cap L_{1,5}=(-2a^2+ac+ad, ac-ad,c-d),\\
$$L_{8,13}: (c-d)x+(2a+c-d)y-2acz=0
\end{array}
$$
and finally
$$P_{19}=L_{8,13}\cap L_{11,14}=(0,4ac^2-4acd,4ac-4ad+2c^2-4cd+2d^2).$$

Almost all points of the construction are triple without any additional conditions. Only $P_2$ and $P_3$ require an additional condition to be triple, namely:
$$4a(a+c-d)=0.$$
As $a\neq 0$, we get $a+c-d=0$. Thus the parametrization space of this configuration is a line and the configuration has a realization over $\mathbb{Q}$. It is
not difficult to check (with help of, e.g., Singular) that the product of all twelve lines (which obviously is in $I_7^{(3)}$) does not belong to $I_7^2$. Thus,
the triple points of this configuration give another rational example of the non-containment of the third symbolic power into the second ordinary power of an
ideal.

For the convenience of the reader, we enclose the Singular script in the Appendix.

\section{Appendix}
To check that the product of all twelve lines of the configuration $C_7$ does not belong to $I_7^2$, and thus $I_7^{(3)}\nsubseteq I_7^2$, the following Singular script may be used.
{\scriptsize
\begin{verbatim}
LIB "elim.lib";
ring R=(32003,a,d),(x,y,z),dp;
option(redSB);
proc rdideal(number p, number q, number r) {
   matrix m[2][3]=p,q,r,x,y,z;
   ideal I=minor(m,2);
   I=std(I);
   return(I);}
proc pline(list P1, list P2) {
   matrix A[3][3]=P1[1],P1[2],P1[3],P2[1],P2[2],P2[3],x,y,z;
   return(det(A));}
ideal P1=rdideal(1,0,0);
ideal P2=rdideal(-1,1,0);
ideal P3=rdideal(1,1,0);
ideal P4=rdideal(0,0,1);
ideal P5=rdideal(a,a,1);
ideal P6=rdideal(-a,a,1);
ideal P7=rdideal(0,2*a,1);
ideal P8=rdideal((d-a),(d-a),1);
ideal P9=rdideal(-2*a+(d-a),(d-a),1);
ideal P10=rdideal(2*a-(d-a),(d-a),1);
ideal P11=rdideal(-(d-a),(d-a),1);
ideal P12=rdideal(d,d,1);
ideal P13=rdideal(-2*a+d,d,1);
ideal P14=rdideal(2*a-d,d,1);
ideal P15=rdideal(-d,d,1);
ideal P16=rdideal(0,4*a*d,4*a-2*(d-a)+2*d);
ideal P17=rdideal(2*(a2)-a*(d-a)-a*d,a*(d-a)-a*d,(d-a)-d);
ideal P18=rdideal(-2*(a2)+a*(d-a)+a*d, a*(d-a)-a*d,(d-a)-d);
ideal P19=rdideal(0,4*a*((d-a)^2)-4*a*(d-a)*d,4*a*(d-a)-4*a*d+2*((d-a)^2)-4*(d-a)*d+2*(d2));
poly pp=(2*(d*z-y))*((d-a)*z-y)*(a*z-y)*(2*a*d*z-2*a*y+(d-a)*x+(d-a)*y-d*x-d*y)*
(2*a*(d-a)*z-2*a*y+(d-a)*x-(d-a)*y-d*x+d*y)*(2*a*z-x-y)*(-x-y)*x*a*(2*a*z+x-y)*
(-2*a*(d-a)*z+2*a*y+(d-a)*x+(d-a)*y-d*x-d*y)*(-2*a*d*z+2*a*y+(d-a)*x-(d-a)*y-d*x+d*y)*(-y+x);
ideal I=intersect(P1,P2,P3,P4,P5,P6,P7,P8,P9,P10,P11,P12,P13,P14,P15,P16,P17,P18,P19);
I=std(I);
reduce(pp,std(I^2));
\end{verbatim}}
\normalsize

\paragraph*{Acknowledgements.}

  We would like to thank warmly Piotr Pokora for drawing our attention to the paper \cite{B-P}. We also thank  Tomasz Szemberg for discussions and remarks.

   The research of Lampa-Baczy\'nska was partially supported by National Science Centre, Poland, grant 2016/23/N/ST1/01363,  the research of Tutaj-Gasi\'nska was partially supported by National Science Centre, Poland, grant
   2014/15/B/ST1/02197.


\bigskip \footnotesize

\bigskip
   \L ucja Farnik,
   Institute of Mathematics, Polish Academy of Sciences, \'Sniadeckich 8,
PL-00-656 Warszawa, Poland
\\
\nopagebreak
   \textit{E-mail address:} \texttt{Lucja.Farnik@gmail.com}

\bigskip
  Jakub Kabat, Magdalena~Lampa-Baczy\'nska,
   Instytut Matematyki UP,
   Podchor\c a\.zych 2,
   PL-30-084 Krak\'ow, Poland
\\
\nopagebreak
  \textit{E-mail address:} \texttt{xxkabat@gmail.com}\\
  \textit{E-mail address:} \texttt{lampa.baczynska@wp.pl}

\bigskip
   Halszka Tutaj-Gasi\'nska,
   Jagiellonian University, Faculty of Mathematics and Computer Scien\-ce, {\L}ojasiewicza~6, PL-30-348 Krak\'ow, Poland
\\
\nopagebreak
   \textit{E-mail address:} \texttt{halszka.tutaj-gasinska@im.uj.edu.pl}


\end{document}